\documentclass[11pt]{article}
\usepackage{graphicx,amsmath,amssymb,mathrsfs}
\usepackage[a4paper,  margin=1in]{geometry}
\usepackage[utf8]{inputenc}
\usepackage[english]{babel}
\usepackage[utf8]{inputenc}
\usepackage{amsmath}
\usepackage{longtable}
\usepackage{amsfonts}
\usepackage{amssymb} 
\usepackage{enumerate}
\usepackage{enumitem} 
\usepackage{authblk}
\usepackage{natbib}
\usepackage[none]{hyphenat}
\usepackage{setspace}
\setcitestyle{authoryear,open={(},close={)}}

\newtheorem{theorem}{Theorem}[section]

\begin{document}

\title{\textbf{Optimal test statistic under normality assumption }}

\author[1]{Nabaneet Das}
\author[2]{Subir K. Bhandari}
\affil[1]{Indian Statistical Institute,Kolkata}
\affil[2]{Indian Statistical Institute,Kolkata}

\maketitle

\begin{abstract}

The idea of an optimal test statistic in the context of simultaneous hypothesis testing was given by \cite{sun2009large} which is the conditional probability of a hypothesis being null given the data. Since we do not have a simplified expression of the statistic, it is impossible to implement the optimal test in more general dependency setup. This note simplifies the expression of optimal test statistic of \cite{sun2009large} under the multivariate normal model. We have considered the model of \cite{xie2011optimal}, where the test statistics are generated from a multivariate normal distribution conditional to the unobserved states of the hypotheses and the states are i.i.d. Bernoulli random variables. While the equivalence of LFDR and optimal test statistic was established under very stringent conditions of \cite{xie2016correction}, the expression obtained in this paper is valid for any covariance matrix and for any fixed $0<p<1$. The optimal procedure is implemented with the help of this expression and the performances have been compared with Benjamini Hochberg method and marginal procedure. 
\end{abstract}

\section{ Introduction }

Dependent observations are frequently encountered in large scale multiple testing problems and they pose a major challenge because of the limitations of the traditional methods which were developed under the assumption of independence. Examples include micro-array experiments where we come across data on thousands of genes and the goal is to separate the 'significant' ones which are very few in number. Analysis of false discovery rate (FDR)  (\cite{benjamini1995controlling}) have been widely used in such cases. Although the original FDR controlling procedure was developed for independent p values, \cite{benjamini2001control} showed that these p-value based procedures are adaptive to certain dependency structures. However, when the proportion of true nulls is relatively small, these procedures often exhibit undesired results (e.g.- too conservative) \\ 
It can be seen that, in dealing with dependent hypotheses, the validity issue has been over emphasized and very few literature are available which actually address the issue of efficiency. \cite{efron2001empirical} introduced local false discovery rate (LFDR) in z-value based testing procedures and studied both size and power (\cite{efron2007size}). \cite{efron2007correlation}, \cite{efron2010correlated} further investigated the effect of correlations on these z-value based procedures and pointed out that, root mean square (rms) of correlations is an important aspect in determining the validity of these z-value based methods. An excellent review of the whole work can be found in \cite{efron2012large}.\\ 
\cite{sun2009large} took a different approach and developed an adaptive multiple testing rule for false discovery control. In their paper, they have used marginal false discovery rate (mFDR) and marginal false non-discovery rate (mFNR) in place of the traditional FDR and FNRs. However, \cite{genovese2002operating} have established that, under the assumption of independence, these are asymptotically the same in the sense that, mFDR = FDR + $O(\frac{1}{ \sqrt{n}})$ and mFNR = FNR +$O(\frac{1}{ \sqrt{n}})$, where $n$ is the number of hypotheses. Such asymptotic equivalence is valid in a more general setting (\cite{xie2011optimal}) of short range dependency structure. \cite{sun2009large} established a one to one correspondence between weighted classification problem and multiple hypothesis testing problem under the monotone ratio condition (MRC) and introduced a new test statistic named local index of significance (LIS) which is optimal in the sense that the test based on this statistic minimizes the mFNR among all methods that control mFDR at a certain level of significance. The optimality of their test statistic is a remarkable development because of the following two reasons.
\begin{enumerate}[label=\Roman*]
    \item  \:   It does not depend on the structure of dependency of the hypotheses.   
    \item  \:   It has been established under MRC condition which is fairly general. As \cite{sun2009large} has highlighted that, the test statistics that are defined on the basis of z-values, such as local false discovery rate (\cite{efron2001empirical}), p-value and the weighted p-value vector (\cite{genovese2006false}) belong to the MRC class. 
\end{enumerate}
The LIS statistic reduces to local false discovery rate (LFDR) under independence and and the optimality of LFDR based procedures of \cite{efron2012large}, \cite{efron2001empirical} is thus established (\cite{sun2007oracle}). However, the closed form expression of this optimal statistic is usually very difficult to find and this poses a major challenge to its application in real data. \cite{sun2009large} considered the hidden Markov model (HMM) where the latent indicator variable of being non-null follows a homogeneous irreducible Markov Chain and developed a recursive method for implementation of the optimal statistic based test. \cite{xie2011optimal} have implemented this test under multivariate normal distribution model. However, their original claim that, the optimal LIS statistic and LFDR is asymptotically the same, only holds under very stringent conditions imposed on model parameters (\cite{xie2016correction}). In this article, we have studied the same model and substantially simplified the test statistic. The reason for considering multivariate Gaussian model is its wide applicability in real life problems and the results proved in this article hold for any positive definite correlation matrix.  

\section{Oracle Decision rule for Multivariate normal model} 
We consider testing n null hypotheses $ H_{01}, ..., H_{0n}$ and for $i=1,2,..,n $\\

$ \theta_i =
\left\{
	\begin{array}{ll}
		1  & \mbox{if } \text{ i-th null hypothesis is false}  \\
	    0  & \mbox{if } \text{Otherwise} 
	\end{array}
\right.$ \\  

Let $ \mathbf{X} = (X_1, ..., X_n) $ be a sequence of test statistics for testing $ \mathbf{H_0} = (H_{01} , \dots , H_{0n}) $.\\ 
In this paper we consider the following model.  
\begin{itemize} 
\item $ \theta_1 , \dots , \theta_n \stackrel{i.i.d.}{\sim} Ber (p) $ for some $ 0 < p < 1 $ 
\item $\mathbf{X} | \mathbf{ \theta} \sim N_n ( k \mathbf{ \theta} , \Sigma ) $ where $ \mathbf{ \theta} = ( \theta_1, \dots, \theta_n ) $ and $ k \neq 0$. Here $ \Sigma $ is any positive definite covariance matrix. 
\end{itemize}  

\subsection{ Discussion on error rate criteria} 

For any multiple testing procedure on these n hypotheses, let V,R,W,A denote the number of false rejections, no. of rejections, no. of false acceptances and no. of acceptances respectively. The false discovery rate (FDR) and marginal false discovery rate (mFDR) are defined as below. 
$$ \text{ FDR } = E \left[ \frac{V}{R} \right]  \: \: \:  \text{      and               }  \: \: \:  mFDR = \frac{ E [ V ] }{ E [ R ] } $$ 

These are versions of type - I error in the context of multiple testing. And, the versions of type -II errors are defined as 

$$ \text{ FNR } = E \left[ \frac{W}{A} \right]  \: \: \:  \text{      and               }  \: \: \:  mFNR = \frac{ E [ W ] }{ E [ A ] } $$

It can be easily shown by Jensen's inequality that, 

$$ \text{ FDR (FNR) } \leq \text{ mFDR (mFNR) } $$
This implies, the methods which aims to control mFDR, tend to be more conservative than the methods controlling FDR. However, \cite{genovese2002operating} has shown that, $ \text{ mFDR (mFNR) } = \text{FDR (FNR) } + O( \frac{1}{ \sqrt{n}} )  $ under independence. \cite{xie2011optimal} have established the asymptotic equivalence of mFDR (mFNR) and FDR (FNR) under short range dependency criterion. 
\begin{theorem}(\cite{xie2011optimal})
Suppose $\mathbf{X} = (X_1, \dots, X_n)$ is a sequence of random variables with same marginal density $ f $ and $ X_i \text{ and } X_j$ are independent if $ |i -j | > n^{ \tau} $ for some $ 0 \leq \tau < 1 $. Let, $ \hat{ \delta}_i = I_{ (S_i \in R } $ be a short-ranged rule to test $H_{i0}$, in the sense that $ S_i $ only depends on the variables that are dependent with $X_i$, $$S_i = S (X_{ i - [n^{ \tau} ] }, ....,X_{ i + [n^{ \tau} ] } ) $$
Further, suppose that, $$P(S_i \in R , \theta_i = 1 ) \geq P(S_i \in R , \theta_i = 0),  $$ and, $$ P(S_i \in R , \theta_i =1 ) > 0 \: \: \forall \: i =1,2,..,n.$$
Then, the FDR(FNR) of the rule $ \hat{ \delta} $ can be approximated by the mFDR(mFNR) in the sense that, $$ mFDR = FDR + O ( \frac{1}{ n^{ 1 - \tau} } ) \: \: \: \text{   and    } \: \: \:mFNR = FNR + O ( \frac{1}{ n^{ 1 - \tau} } ) $$
\end{theorem} 
\textbf{\underline{Note} :-} This asymptotic equivalence does not hold when the correlation matrix is not sparse (e.g. equi-correlated case)\\
In this article, we have considered mFDR and mFNR and derived an optimal test statistic which minimizes mFNR among all methods controlling mFDR at a pre-specified level of significance. 

\subsection{Oracle decision rule for multiple testing problem} 
Consider the weighted classification problem with decision rule $ \delta = ( \delta_1, ..., \delta_n ) \in \{ 0, 1 \}^n$, with $ \delta_i =1 $ if i-th hypothesis is rejected and $ \delta_i = 0 $ otherwise.\\
Consider the loss function 
\begin{equation}\label{1} L_{ \lambda} ( \delta , \theta ) = \frac{1}{n} \sum\limits_{ i=1}^{n} \{ \delta_i (1 - \theta_i) + \lambda \theta_i ( 1 - \delta_i) \}
\end{equation} 
with $ \lambda > 0$ the weight for a false positive result. It is well-known that, if $ g(x | \theta_i = j) $ denotes the density of $x$ when $ \theta_i = j $ ( $ j =0,1 $), then, the classification risk $E[ L_{ \lambda } ( \theta , \delta ) ]  $ is minimized by the Bayes rule $\delta ( \Lambda , \lambda ) = ( \delta_1 , ..., \delta_n ) $, where 

\begin{equation}\label{2}
\delta_i = I \{ \Lambda_i (x ) = \frac{ (1-p) g ( x | \theta_i = 0 ) }{ p g (x | \theta_i =1 ) } < \lambda \}   
\end{equation} 
Alternatively, if the goal is to discover as many significant hypotheses as possible while incurring a relatively low proportion of false positives, we can study a multiple testing problem where the goal is find a decision rule $ \delta $ that has the smallest FNR(mFNR) amoung all FDR(mFDR) procedures at level $\alpha$. \cite{sun2009large}, \cite{xie2011optimal} has shown that, among all procedures controlling mFDR at level $ \alpha$, a procedure which minimizes the mFNR must be of the form $ \delta (\mathbf{T, c} )  = I_{ \mathbf{T } < c \mathbf{1} }= I(T_i < c , \: \: i=1,2,..,n)   $ for some statistic $ \mathbf{T}$ and some real number $c$. (Here $ \mathbf{1}$ denotes the vector with all entries equal to 1) The following theorem explicate the whole idea. 
\begin{theorem}\label{th2.2}( \cite{xie2011optimal} )
Consider the class of decision rules $ \mathscr{D}_s = \{ \delta \: : \: \delta_i = I_{ \Lambda_i < \lambda } , i=1, \dots, n \} $ where $ \mathbf{\Lambda} = \{ \Lambda_1 , \dots , \Lambda_n \}  $ is defined in (\ref{2})  and $ \lambda \in \mathbb{R}$. Given any mFDR level $ \alpha$  and a decision rule 
$$ \delta (S , R ) = \{ I_{ S_1 \in R_1 } , ..., I_{ S_n \in R_n } \} $$ with $ mFDR ( \delta (S, R) ) \leq \alpha $. Then there exists a $ \lambda $ depending on $ \delta (S, R ) $, such that, $\delta(\Lambda, \lambda) \in \mathscr{D}_s$ outperforms $ \delta( S, R) $ in the sense that, 
$$mFDR(\delta ( \Lambda , \lambda ) \leq mFDR (\delta ( S, R) ) \leq \alpha,$$ and 
$$mFNR(\delta ( \Lambda , \lambda ) \leq mFNR (\delta ( S, R) )$$
\end{theorem} 

Theorem \ref{th2.2} implies that, the optimal solution of the multiple testing problem with mFDR and mFNR as the error rate criteria, belongs to the set $\mathscr{D}_s$. Instead of searching for all decision rules, one only needs to search in the collection $\mathscr{D}_s$ for the optimal rule. The following result shows that, for a given $\alpha$, the optimal rule for the multiple testing problem is unique.
\begin{theorem}\label{th2.3} ( \cite{xie2011optimal} )
Consider the optimal decision rule $\delta (\Lambda , \lambda) $ in the weighted classification problem with the loss function (\ref{1}). For any $0<\alpha<1$, there exists a unique $\lambda ( \alpha) $, such that $\delta \{ \Lambda , \lambda ( \alpha)  \}$ controls the mFDR at level $\alpha$ and minimizes the mFNR among all decision rules. 
\end{theorem} 
Theorem \ref{th2.3} gives us the optimal testing rule with mFDR and mFNR as the error rate criteria. It also establishes a one-to-one correspondence between the multiple testing problem and weighted classification problem. However, it is often hard to determine the $\lambda(\alpha) $ corresponding to the given $\alpha$. 
\subsection{Implementation of the optimal test} 
\cite{xie2011optimal} have provided a method to implement the optimal test of theorem 2.2. Define, 
$$T_{OR,i} = P( \theta_i = 0 | x) = \frac{ (1-p) g( x | \theta_i = 0 ) }{ g(x) } $$ 
Clearly, $ T_{OR,i} = \frac{ \Lambda_i}{ 1  + \Lambda_i} $
increases with $\Lambda_i$. Thus, for a given mFDR value $\alpha$, one can rewrite the optimal rule as 
$$\delta_{OR,i}  = \delta ( \Lambda , \lambda ( \alpha) ) = I \left\{ T_{OR, i} < \frac{ \lambda ( \alpha ) }{ 1 + \lambda ( \alpha ) }  \right\}$$ 
Let $T_{OR,(i)} $ denote the i-th order statistic of $T_{OR,i}$ and $H_{0,(i)} $ be the corresponding null hypothesis $(i=1, \dots, n) $ . Then, if $R$ denote the no. of rejections, then 
$$mFDR = E \left[\frac{1}{R} \sum\limits_{i=1}^{R} T_{OR, (i) } \right] $$ 
Then, according to the theorem 5 of \cite{xie2011optimal}, if $p$ and $g$ are known, then the following method controls mFDR at level $\alpha$ : 
\begin{equation}\label{3}
 \text{ Reject all } H_{0, (i) } \text{    for } i=1, \dots, k \: \: \text{    where   } k = \max \left\{ l \: : \: \frac{1}{l} \sum\limits_{i=1}^{l} T_{OR, (i) } \leq \alpha  \right\}   
\end{equation} 
The final oracle rule (\ref{3}) consists of two steps : 
\begin{itemize}
    \item Calculate the oracle statistic $T_{OR,i}$ for $i=1, \dots, n$. 
    \item Rank the statistics and calculate the running averages to determine the cutoff. All hypotheses below the cutoff are rejected. 
\end{itemize}  
However, the major difficulty associated with this optimal test is that the test statistic $T_{OR,i}$ is often very difficult to compute. A simplified expression of this test statistic is hard to find and the model parameters are difficult to estimate under dependent models. In this article, we provide a method for implementation of the optimal test under the multivariate normal model. 

\section{Simplification of the oracle decision rule for multivariate normal model} 
 Under the model specified in section 2, the optimal test statistic can be simplified as follows. 
\begin{theorem}\label{th3.1} 
If $\mathbf{X} \:  | \: \mathbf{\theta} \sim N_n ( k \mathbf{\theta} , \Sigma) $ and $ \theta_1 , \dots , \theta_n \sim^{ i.i.d.} Ber (p) $, then for $ i = 1, \dots, n $  

$$P( \theta_i = 0 \: | \:  \mathbf{X} ) \: = \:  \frac{1}{1 + \frac{pU_i}{1-p}} $$  

Where $ U_i = \exp ( - ( \frac{k^2}{2} t_{i,i}  -  k \sum\limits_{ j=1}^{n} t_{j,i}  x_j )) ( \prod\limits_{j \neq i}  ( pe^{ - k^2 t_{j,i} } + (1-p) )  )  $  and   $\underset{\sim}{t_i} =(t_{1,i} , ..., t_{n,i} ) $ is the $i$-th column of $\Sigma^{-1}$.  
\end{theorem} 
Proof of theorem \ref{th3.1} is given in appendix. \\ 

\textbf{\underline{Remarks}} :- 
 While the equivalence of the joint conditional probability of \cite{xie2011optimal} and \cite{xie2016correction} was established under very restrictive assumptions on the correlation matrix, theorem \ref{th3.1} sufficiently simplifies the conditional probability for any covariance matrix. The result of \ref{th3.1} enables us to implement the Oracle decision rule for any $p, \Sigma$. We have performed extensive simulations with different combinations of $p$ and $\Sigma $ and compared the observed value of FDR and FNRs some of which will be discussed here. 

\subsection{Simulation Studies}  
In this section, we evaluate the performance of the oracle rule and compare with the BH procedure and the marginal procedure mentioned in \cite{xie2011optimal}. We have evaluated the empirical FDR, FNR and also the number of rejections. In our simulations, we assumed a multivariate normal model : 
$$ X | \theta \sim N( c \theta , \Sigma) $$ 
where $\theta_i$ follows Bernoulli($p$). Under this model, the non-null distribution has mean $c$ and $\Sigma$ is a correlation matrix. For our simulations, we have considered $c= 2.5$ and $\alpha = 0.05$. Our objective is to assess the performance of these methods when there is sufficient deviation from independence. In the first case, we have considered equicorrelated $\Sigma$. In all the simulations, number of hypotheses ($n$) have been considered to be 5000 and they are run on 10 combinations of proportion of non-null $p = 0.01, 0.02, ...., 0.1$. In order to assess the performance under sufficient deviation from independence, seven cominations of correlation have been considered ($\rho = 0.2,0.3,0.4,0.5,0.6,0.7,0.8$). The results suggest that, the Oracle procedure is least conservative among the three procedures in terms of FDR. It is interesting to note that, the FDR of the Oracle rule always lies within the prescribed limit of 0.05. Maintaining this upper bound on FDR, there is a substantial gain in the FNR over both BH and marginal procedure. It is interesting to note that, the marginal procedure becomes more conservative than the other two methods. However, the FNR of the marginal procedure remains similar to the BH procedure which suggests a possibility of improvement of this method and that is achieved by considering the information of joint distribution in the Oracle procedure. The conservative nature of the marginal procedure in comparison to the BH method is possibly due to the objective of controlling mFDR instead of FDR. Since the FNR remains equivalent for these two methods, careful examination of the class $\mathscr{D}_s$ may provide a significantly better test statistic.\\ 
\begin{center}{}
\begin{longtable}{|l|l|l|l|l|}
\caption[FDRs of the three methods]{FDRs of the three methods}  \\

\hline \multicolumn{1}{|c|}{\textbf{Marginal Procedure}} & \multicolumn{1}{c|}{\textbf{BH procedure}} & 
\multicolumn{1}{c|}{\textbf{Oracle procedure}} & 
\multicolumn{1}{c|}{\textbf{p}} & 
\multicolumn{1}{c|}{\textbf{Correlation}}\\ \hline 
\endfirsthead

\multicolumn{5}{c}%
{{\bfseries \tablename\ \thetable{} (FDRs)-- continued from previous page}} \\
\hline \multicolumn{1}{|c|}{\textbf{Marginal Procedure}} & \multicolumn{1}{c|}{\textbf{BH procedure}} & 
\multicolumn{1}{c|}{\textbf{Oracle procedure}} & 
\multicolumn{1}{c|}{\textbf{p}} & 
\multicolumn{1}{c|}{\textbf{Correlation}} \\ \hline 
\endhead

\hline \multicolumn{5}{|r|}{{Continued on next page}} \\ \hline
\endfoot

\hline \hline
\endlastfoot

0.016177118        & 0.042388341  & 0.042601571      & 0.01 & 0.2         \\
0.010775516        & 0.037959313  & 0.044011605      & 0.01 & 0.3         \\
0.007105661        & 0.032570762  & 0.044454383      & 0.01 & 0.4         \\
0.004964538        & 0.028604701  & 0.043996202      & 0.01 & 0.5         \\
0.003348643        & 0.022948649  & 0.043062505      & 0.01 & 0.6         \\
0.002811004        & 0.019151864  & 0.041333631      & 0.01 & 0.7         \\
0.002425759        & 0.017037651  & 0.038867168      & 0.01 & 0.8         \\
\hline
0.022105103        & 0.042394141  & 0.044522453      & 0.02 & 0.2         \\
0.015325569        & 0.039289712  & 0.044120863      & 0.02 & 0.3         \\
0.010682572        & 0.034905345  & 0.043407766      & 0.02 & 0.4         \\
0.007820876        & 0.030941068  & 0.042582097      & 0.02 & 0.5         \\
0.006264827        & 0.026842743  & 0.041341872      & 0.02 & 0.6         \\
0.005166191        & 0.021776192  & 0.039653903      & 0.02 & 0.7         \\
0.004997578        & 0.018555037  & 0.040080633      & 0.02 & 0.8         \\
\hline
0.02472268         & 0.04238432   & 0.043038606      & 0.03 & 0.2         \\
0.018112997        & 0.039890416  & 0.042364845      & 0.03 & 0.3         \\
0.013348423        & 0.036379311  & 0.041606332      & 0.03 & 0.4         \\
0.010280679        & 0.032402211  & 0.040305484      & 0.03 & 0.5         \\
0.008641967        & 0.028717477  & 0.038804942      & 0.03 & 0.6         \\
0.007632851        & 0.024204622  & 0.037183505      & 0.03 & 0.7         \\
0.00707045         & 0.019213027  & 0.039281896      & 0.03 & 0.8         \\
\hline
0.026978791        & 0.042743758  & 0.04154021       & 0.04 & 0.2         \\
0.020255192        & 0.040077515  & 0.040665379      & 0.04 & 0.3         \\
0.015368755        & 0.036617453  & 0.039435394      & 0.04 & 0.4         \\
0.012312905        & 0.033511309  & 0.037921364      & 0.04 & 0.5         \\
0.01086395         & 0.030662412  & 0.036098383      & 0.04 & 0.6         \\
0.009645082        & 0.025665106  & 0.03433813       & 0.04 & 0.7         \\
0.009102326        & 0.020413258  & 0.038163127      & 0.04 & 0.8         \\
\hline
0.028320129        & 0.042371924  & 0.039951502      & 0.05 & 0.2         \\
0.021869967        & 0.040036222  & 0.038849095      & 0.05 & 0.3         \\
0.017610836        & 0.038088187  & 0.037465988      & 0.05 & 0.4         \\
0.014300737        & 0.034324064  & 0.035646333      & 0.05 & 0.5         \\
0.012634147        & 0.031103329  & 0.033516813      & 0.05 & 0.6         \\
0.011670306        & 0.026734704  & 0.031802998      & 0.05 & 0.7         \\
0.011271098        & 0.022000566  & 0.036869315      & 0.05 & 0.8         \\
\hline
0.029516371        & 0.042079613  & 0.038255191      & 0.06 & 0.2         \\
0.023234026        & 0.039874208  & 0.036971488      & 0.06 & 0.3         \\
0.018900139        & 0.037576243  & 0.035504738      & 0.06 & 0.4         \\
0.016149836        & 0.035093506  & 0.033448897      & 0.06 & 0.5         \\
0.014147447        & 0.03138494   & 0.031066482      & 0.06 & 0.6         \\
0.013427973        & 0.027561753  & 0.029292927      & 0.06 & 0.7         \\
0.013240909        & 0.023446602  & 0.035595795      & 0.06 & 0.8         \\
\hline
0.030703506        & 0.041852071  & 0.036775754      & 0.07 & 0.2         \\
0.024776596        & 0.040187421  & 0.035281783      & 0.07 & 0.3         \\
0.02045121         & 0.037793241  & 0.033618521      & 0.07 & 0.4         \\
0.017519457        & 0.035023497  & 0.031401037      & 0.07 & 0.5         \\
0.015842256        & 0.032099682  & 0.028824389      & 0.07 & 0.6         \\
0.015019501        & 0.028258038  & 0.026953824      & 0.07 & 0.7         \\
0.014502993        & 0.023252437  & 0.034348857      & 0.07 & 0.8         \\
\hline
0.031997691        & 0.04225768   & 0.035212564      & 0.08 & 0.2         \\
0.025856448        & 0.039742567  & 0.033650868      & 0.08 & 0.3         \\
0.021890963        & 0.038079083  & 0.031651006      & 0.08 & 0.4         \\
0.019276162        & 0.035837497  & 0.029363319      & 0.08 & 0.5         \\
0.017228129        & 0.032270158  & 0.026718071      & 0.08 & 0.6         \\
0.017065663        & 0.029549042  & 0.024721255      & 0.08 & 0.7         \\
0.016453646        & 0.024633957  & 0.033020328      & 0.08 & 0.8         \\
\hline
0.032438114        & 0.041155868  & 0.033735288      & 0.09 & 0.2         \\
0.026846327        & 0.039445026  & 0.032056356      & 0.09 & 0.3         \\
0.02284317         & 0.037435289  & 0.030007486      & 0.09 & 0.4         \\
0.020452384        & 0.035578941  & 0.027485875      & 0.09 & 0.5         \\
0.019038378        & 0.032965301  & 0.024638386      & 0.09 & 0.6         \\
0.018614293        & 0.029862288  & 0.022581026      & 0.09 & 0.7         \\
0.017985636        & 0.025240484  & 0.031805517      & 0.09 & 0.8         \\
\hline
0.03336696         & 0.041015676  & 0.032376719      & 0.1  & 0.2         \\
0.028046145        & 0.039468922  & 0.03054991       & 0.1  & 0.3         \\
0.024476314        & 0.038105886  & 0.028347872      & 0.1  & 0.4         \\
0.021867929        & 0.03587906   & 0.025697372      & 0.1  & 0.5         \\
0.019978473        & 0.032617402  & 0.022729948      & 0.1  & 0.6         \\
0.019614439        & 0.029609253  & 0.020594574      & 0.1  & 0.7         \\
0.019335304        & 0.025153741  & 0.03046926       & 0.1  & 0.8        \\
\end{longtable}
\end{center}
From the FDRs, it is clear that, 
\begin{itemize}
    \item All the three methods (especially BH procedure) become more and more conservative with increasing value of the correlation. 
    \item Marginal procedure tend to be the most conservative among the other three methods for higher correlations. (i.e. $\text{FDR}_{MP} \leq \text{ FDR}_{BH} \leq \text{FDR}_{OP}$. However, slight exceptions can be observed for smaller correlations and higher $p$ where BH procedure has slightly higher FDR among all. 
\end{itemize}
With the above mentioned observations on FDR, it is imperative to note the FNRs of these three methods. 
\begin{center}
\begin{longtable}{|l|l|l|l|l|}
\caption[FNRs of the three methods]{FNRs of the three methods}  \\

\hline \multicolumn{1}{|c|}{\textbf{Marginal Procedure}} & \multicolumn{1}{c|}{\textbf{BH procedure}} & 
\multicolumn{1}{c|}{\textbf{Oracle procedure}} & 
\multicolumn{1}{c|}{\textbf{p}} & 
\multicolumn{1}{c|}{\textbf{Correlation}}\\ \hline 
\endfirsthead

\multicolumn{5}{c}%
{{\bfseries \tablename\ \thetable{}(FNRs) -- continued from previous page}} \\
\hline \multicolumn{1}{|c|}{\textbf{Marginal Procedure}} & \multicolumn{1}{c|}{\textbf{BH procedure}} & 
\multicolumn{1}{c|}{\textbf{Oracle procedure}} & 
\multicolumn{1}{c|}{\textbf{p}} & 
\multicolumn{1}{c|}{\textbf{Correlation}} \\ \hline 
\endhead

\hline \multicolumn{5}{|r|}{{Continued on next page}} \\ \hline
\endfoot

\hline \hline
\endlastfoot

0.009384108        & 0.008939474  & 0.008485681      & 0.01 & 0.2         \\
0.00938644         & 0.008888759  & 0.007651516      & 0.01 & 0.3         \\
0.009378513        & 0.008838835  & 0.006468921      & 0.01 & 0.4         \\
0.009363355        & 0.008790743  & 0.004887012      & 0.01 & 0.5         \\
0.009386034        & 0.008774157  & 0.003013366      & 0.01 & 0.6         \\
0.009381393        & 0.008718203  & 0.001196976      & 0.01 & 0.7         \\
0.009386836        & 0.008648042  & 0.000144564      & 0.01 & 0.8         \\
\hline
0.017667805        & 0.017049137  & 0.01537233       & 0.02 & 0.2         \\
0.017673984        & 0.016923198  & 0.013481906      & 0.02 & 0.3         \\
0.017701522        & 0.016848171  & 0.011025832      & 0.02 & 0.4         \\
0.017726805        & 0.016764811  & 0.008004701      & 0.02 & 0.5         \\
0.017709012        & 0.016620481  & 0.004660291      & 0.02 & 0.6         \\
0.017780044        & 0.016541724  & 0.001697489      & 0.02 & 0.7         \\
0.017855209        & 0.016391572  & 0.000162737      & 0.02 & 0.8         \\
\hline
0.025252815        & 0.024618774  & 0.021579133      & 0.03 & 0.2         \\
0.025285156        & 0.024444385  & 0.018668056      & 0.03 & 0.3         \\
0.025307303        & 0.024297564  & 0.014977791      & 0.03 & 0.4         \\
0.025368444        & 0.024186403  & 0.010641615      & 0.03 & 0.5         \\
0.025417939        & 0.024040624  & 0.006023362      & 0.03 & 0.6         \\
0.025528645        & 0.023880292  & 0.002089584      & 0.03 & 0.7         \\
0.025651744        & 0.023590965  & 0.00017367       & 0.03 & 0.8         \\
\hline
0.032238775        & 0.031678249  & 0.027352014      & 0.04 & 0.2         \\
0.032300073        & 0.031491575  & 0.023427774      & 0.04 & 0.3         \\
0.03240079         & 0.031374052  & 0.018616471      & 0.04 & 0.4         \\
0.032477425        & 0.031201011  & 0.013049375      & 0.04 & 0.5         \\
0.032534081        & 0.03095588   & 0.007271547      & 0.04 & 0.6         \\
0.032704843        & 0.030741584  & 0.002446214      & 0.04 & 0.7         \\
0.03309348         & 0.030535566  & 0.000182017      & 0.04 & 0.8         \\
\hline
0.038870145        & 0.038483254  & 0.032859021      & 0.05 & 0.2         \\
0.038908465        & 0.038235336  & 0.027954742      & 0.05 & 0.3         \\
0.039014219        & 0.038067704  & 0.022069251      & 0.05 & 0.4         \\
0.039193163        & 0.037966031  & 0.015342851      & 0.05 & 0.5         \\
0.039259985        & 0.037618644  & 0.008479741      & 0.05 & 0.6         \\
0.039595532        & 0.037439551  & 0.002798795      & 0.05 & 0.7         \\
0.040013936        & 0.037044964  & 0.000189161      & 0.05 & 0.8         \\
\hline
0.04510414         & 0.044968515  & 0.038171907      & 0.06 & 0.2         \\
0.045211559        & 0.044774224  & 0.03234459       & 0.06 & 0.3         \\
0.045333742        & 0.044582031  & 0.025414657      & 0.06 & 0.4         \\
0.045462535        & 0.044358697  & 0.01758729       & 0.06 & 0.5         \\
0.045639445        & 0.044021015  & 0.009663691      & 0.06 & 0.6         \\
0.045999712        & 0.043735292  & 0.003167917      & 0.06 & 0.7         \\
0.046636465        & 0.043322141  & 0.000197412      & 0.06 & 0.8         \\
\hline
0.050986433        & 0.05114727   & 0.043355167      & 0.07 & 0.2         \\
0.051128908        & 0.050993789  & 0.036619159      & 0.07 & 0.3         \\
0.05123691         & 0.050737724  & 0.028681921      & 0.07 & 0.4         \\
0.051464837        & 0.050539536  & 0.019833372      & 0.07 & 0.5         \\
0.051749431        & 0.050272812  & 0.010874287      & 0.07 & 0.6         \\
0.052143392        & 0.049858717  & 0.003531483      & 0.07 & 0.7         \\
0.052947299        & 0.049431048  & 0.000204346      & 0.07 & 0.8         \\
\hline
0.056536873        & 0.057057084  & 0.048426578      & 0.08 & 0.2         \\
0.056758056        & 0.056963728  & 0.040844216      & 0.08 & 0.3         \\
0.056998609        & 0.056851589  & 0.031952433      & 0.08 & 0.4         \\
0.057243247        & 0.05661948   & 0.022079845      & 0.08 & 0.5         \\
0.057478277        & 0.05617847   & 0.012100794      & 0.08 & 0.6         \\
0.058067235        & 0.055854904  & 0.003927636      & 0.08 & 0.7         \\
0.059003171        & 0.055261631  & 0.000213031      & 0.08 & 0.8         \\
\hline
0.062053951        & 0.063041149  & 0.053491892      & 0.09 & 0.2         \\
0.062156781        & 0.062813482  & 0.045037258      & 0.09 & 0.3         \\
0.062395084        & 0.062622997  & 0.035225876      & 0.09 & 0.4         \\
0.062544115        & 0.06224676   & 0.024371961      & 0.09 & 0.5         \\
0.063065465        & 0.062044864  & 0.013404678      & 0.09 & 0.6         \\
0.063559517        & 0.061412199  & 0.00435216       & 0.09 & 0.7         \\
0.064711954        & 0.060816922  & 0.000219888      & 0.09 & 0.8         \\
\hline
0.067187996        & 0.068644091  & 0.058500897      & 0.1  & 0.2         \\
0.06733126         & 0.068470557  & 0.049254304      & 0.1  & 0.3         \\
0.067512376        & 0.068229779  & 0.038542111      & 0.1  & 0.4         \\
0.067828866        & 0.067951705  & 0.026689565      & 0.1  & 0.5         \\
0.068306286        & 0.067548909  & 0.014756385      & 0.1  & 0.6         \\
0.069211136        & 0.067231416  & 0.004816067      & 0.1  & 0.7         \\
0.070454012        & 0.066548327  & 0.00022843       & 0.1  & 0.8         \\
\end{longtable}
\end{center} 
As per the optimality of Oracle procedure, it has the lowest FNR among all. It is interesting to note that, while the marginal procedure was the most conservative in terms of FDR, its FNR is nearly equivalent (or even better in some cases) to the BH procedure. It is again reminded that, we are controlling the mFDR(mFNR) instead of FDR(FNR). The results suggest that there is a scope of further improvement in the class $ \mathscr{D}_s = \{ \delta \: : \: \delta_i = I_{ \Lambda_i < \lambda } , i=1, \dots, n \} $ if $ \Lambda $ and $ \lambda $ can be chosen properly.\\
Examining the FDRs and FNRs does not entirely describe how conservative a method is. We know that, these methods become conservative with increasing value of correlation. To examine this, we have also tabulated the no. of rejections of these three methods in different combinations of $p$ and correlation. 
\begin{center}
\begin{longtable}{|l|l|l|l|l|}
\caption[No. of rejections of the three methods]{No. of rejections of the three methods}  \\

\hline \multicolumn{1}{|c|}{\textbf{Marginal Procedure}} & \multicolumn{1}{c|}{\textbf{BH procedure}} & 
\multicolumn{1}{c|}{\textbf{Oracle procedure}} & 
\multicolumn{1}{c|}{\textbf{p}} & 
\multicolumn{1}{c|}{\textbf{Correlation}}\\ \hline 
\endfirsthead

\multicolumn{5}{c}%
{{\bfseries \tablename\ \thetable{}(No. of rejections) -- continued from previous page}} \\
\hline \multicolumn{1}{|c|}{\textbf{Marginal Procedure}} & \multicolumn{1}{c|}{\textbf{BH procedure}} & 
\multicolumn{1}{c|}{\textbf{Oracle procedure}} & 
\multicolumn{1}{c|}{\textbf{p}} & 
\multicolumn{1}{c|}{\textbf{Correlation}} \\ \hline 
\endhead

\hline \multicolumn{5}{|r|}{{Continued on next page}} \\ \hline
\endfoot

\hline \hline
\endlastfoot

3                  & 8            & 8                & 0.01 & 0.2         \\
3                  & 13           & 12               & 0.01 & 0.3         \\
3                  & 19           & 19               & 0.01 & 0.4         \\
3                  & 28           & 27               & 0.01 & 0.5         \\
3                  & 36           & 37               & 0.01 & 0.6         \\
3                  & 52           & 46               & 0.01 & 0.7         \\
3                  & 70           & 51               & 0.01 & 0.8         \\
\hline
13                 & 19           & 25               & 0.02 & 0.2         \\
13                 & 25           & 34               & 0.02 & 0.3         \\
12                 & 31           & 47               & 0.02 & 0.4         \\
12                 & 40           & 63               & 0.02 & 0.5         \\
12                 & 51           & 80               & 0.02 & 0.6         \\
12                 & 62           & 95               & 0.02 & 0.7         \\
12                 & 83           & 103              & 0.02 & 0.8         \\
\hline
26                 & 34           & 45               & 0.03 & 0.2         \\
26                 & 40           & 60               & 0.03 & 0.3         \\
26                 & 47           & 80               & 0.03 & 0.4         \\
25                 & 55           & 102              & 0.03 & 0.5         \\
25                 & 67           & 125              & 0.03 & 0.6         \\
25                 & 79           & 145              & 0.03 & 0.7         \\
24                 & 95           & 155              & 0.03 & 0.8         \\
\hline
42                 & 51           & 68               & 0.04 & 0.2         \\
42                 & 57           & 89               & 0.04 & 0.3         \\
42                 & 64           & 113              & 0.04 & 0.4         \\
41                 & 72           & 142              & 0.04 & 0.5         \\
41                 & 85           & 171              & 0.04 & 0.6         \\
40                 & 96           & 195              & 0.04 & 0.7         \\
38                 & 111          & 207              & 0.04 & 0.8         \\
\hline
61                 & 69           & 92               & 0.05 & 0.2         \\
61                 & 75           & 118              & 0.05 & 0.3         \\
61                 & 83           & 149              & 0.05 & 0.4         \\
60                 & 92           & 183              & 0.05 & 0.5         \\
60                 & 104          & 217              & 0.05 & 0.6         \\
58                 & 117          & 244              & 0.05 & 0.7         \\
56                 & 130          & 259              & 0.05 & 0.8         \\
\hline
83                 & 89           & 118              & 0.06 & 0.2         \\
82                 & 95           & 149              & 0.06 & 0.3         \\
82                 & 103          & 184              & 0.06 & 0.4         \\
81                 & 113          & 223              & 0.06 & 0.5         \\
80                 & 123          & 262              & 0.06 & 0.6         \\
78                 & 136          & 294              & 0.06 & 0.7         \\
75                 & 153          & 310              & 0.06 & 0.8         \\
\hline
106                & 111          & 145              & 0.07 & 0.2         \\
105                & 118          & 180              & 0.07 & 0.3         \\
105                & 125          & 220              & 0.07 & 0.4         \\
104                & 134          & 264              & 0.07 & 0.5         \\
102                & 145          & 308              & 0.07 & 0.6         \\
100                & 157          & 343              & 0.07 & 0.7         \\
96                 & 170          & 361              & 0.07 & 0.8         \\
\hline
131                & 135          & 172              & 0.08 & 0.2         \\
130                & 141          & 212              & 0.08 & 0.3         \\
129                & 149          & 257              & 0.08 & 0.4         \\
128                & 158          & 305              & 0.08 & 0.5         \\
126                & 166          & 353              & 0.08 & 0.6         \\
124                & 182          & 392              & 0.08 & 0.7         \\
119                & 195          & 413              & 0.08 & 0.8         \\
\hline
157                & 158          & 200              & 0.09 & 0.2         \\
157                & 165          & 244              & 0.09 & 0.3         \\
155                & 172          & 293              & 0.09 & 0.4         \\
155                & 182          & 346              & 0.09 & 0.5         \\
152                & 192          & 398              & 0.09 & 0.6         \\
150                & 207          & 440              & 0.09 & 0.7         \\
143                & 219          & 464              & 0.09 & 0.8         \\
\hline
186                & 184          & 228              & 0.1  & 0.2         \\
185                & 190          & 276              & 0.1  & 0.3         \\
184                & 198          & 329              & 0.1  & 0.4         \\
182                & 206          & 387              & 0.1  & 0.5         \\
179                & 215          & 443              & 0.1  & 0.6         \\
175                & 227          & 488              & 0.1  & 0.7         \\
167                & 241          & 515              & 0.1  & 0.8          \\
\end{longtable}
\end{center}
No. of rejections for the Oracle procedure is significantly higher than the other three methods and hence this is the least conservative among all. However, the equicorrelated $\Sigma$ is an unlikely scenario in real life applications. We only considered this in order to generate a scenario which is substantially different from the independent setup and compare the performances of the methods. Now we present the results on block diagonal correlation matrix. Here we have divided the correlation matrix in four blocks of equicorrelated matrices with correlation $0.15, 0.25, 0.5, 0.75$. The results again suggest that, the Oracle Procedure is least conservative among the three methods in terms of FDR while maintaining the prescribed limit of 0.05. No. of rejections for Oracle procedure is significantly higher than the other two and the gain in power is also noteworthy.

\begin{center}
\begin{table}[ht]
\begin{tabular}{|l|l|l|l|l|l|l|l|}
\hline
Marginal Procedure & BH Procedure & Oracle Procedure & p & $\rho_1$ & $\rho_2$ & $\rho_3$ & $\rho_4$\\
\hline
0.018467           & 0.036869     & 0.043467         & 0.01 & 0.25   & 0.5    & 0.15   & 0.75   \\
0.027403           & 0.038453     & 0.04318          & 0.02 & 0.25   & 0.5    & 0.15   & 0.75   \\
0.030884           & 0.039435     & 0.041572         & 0.03 & 0.25   & 0.5    & 0.15   & 0.75   \\
0.032916           & 0.039906     & 0.039776         & 0.04 & 0.25   & 0.5    & 0.15   & 0.75   \\
0.034597           & 0.040304     & 0.037964         & 0.05 & 0.25   & 0.5    & 0.15   & 0.75   \\
0.035802           & 0.04027      & 0.03615          & 0.06 & 0.25   & 0.5    & 0.15   & 0.75   \\
0.037047           & 0.040448     & 0.034459         & 0.07 & 0.25   & 0.5    & 0.15   & 0.75   \\
0.038028           & 0.040457     & 0.032838         & 0.08 & 0.25   & 0.5    & 0.15   & 0.75   \\
0.038755           & 0.040247     & 0.031275         & 0.09 & 0.25   & 0.5    & 0.15   & 0.75   \\
0.039511           & 0.040099     & 0.029757         & 0.1  & 0.25   & 0.5    & 0.15   & 0.75   \\ 
\hline

\end{tabular}
\caption{FDRs in block diagonal case}
\label{FDR_1} 
\end{table}
\end{center}

\vspace{0.1 in}

\begin{center}
\begin{table}[ht]
\begin{tabular}{|l|l|l|l|l|l|l|l|}
\hline
Marginal Procedure & BH Procedure & Oracle Procedure & p & $\rho_1$ & $\rho_2$ & $\rho_3$ & $\rho_4$\\
\hline
0.00941            & 0.009047     & 0.005234         & 0.01 & 0.25   & 0.5    & 0.15   & 0.75   \\
0.017709           & 0.017274     & 0.009103         & 0.02 & 0.25   & 0.5    & 0.15   & 0.75   \\
0.025293           & 0.024921     & 0.012522         & 0.03 & 0.25   & 0.5    & 0.15   & 0.75   \\
0.032339           & 0.032104     & 0.015691         & 0.04 & 0.25   & 0.5    & 0.15   & 0.75   \\
0.038957           & 0.038926     & 0.018713         & 0.05 & 0.25   & 0.5    & 0.15   & 0.75   \\
0.045229           & 0.045462     & 0.021645         & 0.06 & 0.25   & 0.5    & 0.15   & 0.75   \\
0.051179           & 0.051733     & 0.024527         & 0.07 & 0.25   & 0.5    & 0.15   & 0.75   \\
0.056849           & 0.057763     & 0.027396         & 0.08 & 0.25   & 0.5    & 0.15   & 0.75   \\
0.062316           & 0.063647     & 0.030275         & 0.09 & 0.25   & 0.5    & 0.15   & 0.75   \\
0.067539           & 0.069341     & 0.033173         & 0.1  & 0.25   & 0.5    & 0.15   & 0.75   \\
\hline

\end{tabular}
\caption{FNRs in block diagonal case}
\label{FNR_1} 
\end{table}
\end{center}

\vspace{0.1 in} 

\begin{center}
\begin{table}[ht]
\begin{tabular}{|l|l|l|l|l|l|l|l|}
\hline
Marginal Procedure & BH Procedure & Oracle Procedure & p & $\rho_1$ & $\rho_2$ & $\rho_3$ & $\rho_4$\\
\hline
3                  & 9            & 25               & 0.01 & 0.25   & 0.5    & 0.15   & 0.75   \\
12                 & 19           & 57               & 0.02 & 0.25   & 0.5    & 0.15   & 0.75   \\
26                 & 33           & 92               & 0.03 & 0.25   & 0.5    & 0.15   & 0.75   \\
42                 & 49           & 129              & 0.04 & 0.25   & 0.5    & 0.15   & 0.75   \\
61                 & 67           & 166              & 0.05 & 0.25   & 0.5    & 0.15   & 0.75   \\
82                 & 87           & 204              & 0.06 & 0.25   & 0.5    & 0.15   & 0.75   \\
106                & 108          & 242              & 0.07 & 0.25   & 0.5    & 0.15   & 0.75   \\
131                & 131          & 280              & 0.08 & 0.25   & 0.5    & 0.15   & 0.75   \\
157                & 156          & 318              & 0.09 & 0.25   & 0.5    & 0.15   & 0.75   \\
186                & 181          & 357              & 0.1  & 0.25   & 0.5    & 0.15   & 0.75   \\
\hline

\end{tabular}
\caption{No. of rejections in block diagonal case}
\label{N_Rej_1} 
\end{table}
\end{center} 

\vspace{0.1 in}

As mentioned earlier, Oracle procedure exploits the information of joint distribution unlike the marginal and BH procedure. The results from the simulation studies have shown a significant improvement in FNR and the no. of rejections in exchange of very little sacrifice in FDR. Hence, it is interesting to explore the class $\mathscr{D} = \{  I_{ X_i > c } \: i =1, \dots , n \}$ and to search for a different choice of $c$ which can provide further improvement. Also, implementation of the optimal procedure under a more general dependency setup (e.g. m-dependent structure) is still a challenging open problem. 

\section{Annexure}
\subsection{Proof of theorem \ref{th3.1}}
Let $ f(\mathbf{x}, \theta ) $ denote the value of $N( k \theta , \Sigma ) $ density at $ \mathbf{x}$.  Then, 
$$ P( \theta_i = 0 | \mathbf{x} ) = \frac{ (1-p) E_{ \theta_{0,i} } [ f  ( \mathbf{x} , \theta_{0,i}) ]}{ E_{ \theta } [ f  ( \mathbf{x} , \theta) ] } 
 $$  
 Where $ \theta_{0,i} $ has $0$ in it's $i$-th place.  \\ 
 Observe that,  
 $$
 \frac{ (1-p) E_{ \theta_{0,i} } [ f  ( \mathbf{x} , \theta_{0,i}) ]}{ E_{ \theta } [ f  ( \mathbf{x} , \theta) ] } \\
  = \frac{ (1-p) \sum\limits_{ \theta_{0,i} }  f  ( \mathbf{x} , \theta_{0,i}) h_{0} ( \theta_{0,i} ) }{ 
  (1-p) \sum\limits_{ \theta_{0,i} }  f  ( \mathbf{x} , \theta_{0,i}) h_{0} (\theta_{0,i} ) + p  \sum\limits_{ \theta_{1,i} }  f  ( \mathbf{x} , \theta_{1,i}) h_{1} ( \theta_{1,i} ) }  
 $$ 
 Where $ \theta_{1,i} $ has $1$ in it's $i$-th place and $h_0$ and $h_1$ are the joint p.m.f.s of $\theta $ given $\theta_i =0 $ and $ \theta_i = 1 $ respectively. In particular, $ \theta_{1,i} = \theta_{0,i} +  e_i   $  ($e_i$ is the vector with 1 in the i-th place and 0 elsewhere) \\ 
  Let $ B = \sum\limits_{ \theta_{1,i} }  f  ( \mathbf{x} , \theta_{1,i}) h_{1} ( \theta_{1,i} ) $  and $A = \sum\limits_{ \theta_{0,i} }  f  ( \mathbf{x} , \theta_{0,i}) h_{0} (\theta_{0,i} ) $. \\

 Then, $  P( \theta_i = 0 | \mathbf{x} ) = \frac{1}{ 1 + \frac{ pB}{(1-p)A}  }  = $  (A monotone function in $ \frac{B}{A}$)  \\ 
 Since $ \theta_i$'s are i.i.d., we must have $h_0 = h_1$. Putting $ \theta_{1,i}  =  \theta_{0,i} + e_i $, we get, 
 
$$ f ( \mathbf{x} , \theta_{1,i})  = f(x , \theta_{0,i}) \exp( - \frac{k^2}{2} e_{i}^T \Sigma^{-1} e_i ) \exp (  k (\mathbf{x} - k \theta_{0,i})^T \Sigma^{-1} e_i) $$

Let  $\underset{\sim}{t_i} =(t_{1,i} , ..., t_{n,i} ) $ be the $i$-th column of $\Sigma^{-1}$. Then, $e_{i}^T \Sigma^{-1} e_i = t_{i,i}$ and $ \mathbf{x}^T \Sigma^{-1} e_i = \sum\limits_{ j =1}^{n} t_{ j,i} x_j $. This implies, 
$$ f  \mathbf{x}  , \theta_{1,i})  = f( \mathbf{x} , \theta_{0,i}) \exp( - (\frac{k^2}{2} t_{i,i}  - k \sum\limits_{j=1}^{n} t_{j,i} x_j) ) (\exp (-k^2 \theta_{0,i}^T \Sigma^{-1} e_i))$$ 
Observe that, $$ \frac{\sum\limits_{ \theta_{1,i} }  f  ( \mathbf{x} , \theta_{1,i}) h_{1} ( \theta_{1,i} )}{\sum\limits_{ \theta_{0,i} }  f  ( \mathbf{x} , \theta_{0,i}) h_{0} (\theta_{0,i} )} = \frac{ E_{ \theta_{1,i} } [ f( \mathbf{x}, \theta_{1,i} ) ]}{ E_{ \theta_{0,i} } [ f ( \mathbf{x}, \theta_{0,i}) ]} =  \exp( - (\frac{k^2}{2} t_{i,i}  - k \sum\limits_{j=1}^{n} t_{j,i} x_j) )\frac{ E_{ \theta_{ 0,i}} [f( \mathbf{x} , \theta_{0,i}) \exp ( -k^2 \theta_{0,i}^T \Sigma^{-1} e_i)]}{E_{ \theta_{0,i} } [ f ( \mathbf{x}, \theta_{0,i}) ]} $$ 
Note that, $ f(\mathbf{x}, \theta_{0,i} ) = g ( \mathbf{x} - \theta_{ 0,i} ) $ = A function of $( \mathbf{x } - k\theta_{0,i})$. 
As per our model $( \mathbf{X} - \theta_{0,i})$ is independent of $ \theta_{0,i}$ and hence, we can say that, 
$$ E_{ \theta_{ 0,i}} [f( \mathbf{x} , \theta_{0,i}) \exp ( - k^2  \theta_{0,i}^T \Sigma^{-1} e_i)]
= E_{ \theta_{ 0,i}} [f( \mathbf{x} , \theta_{0,i}) ] C $$
And thus, 
$$ \frac{B}{A} = \exp( - (\frac{k^2}{2} t_{i,i}  - k \sum\limits_{j=1}^{n} t_{j,i} x_j) )  E_{ \theta_{0,i}} [ \exp ( - k^2 \theta_{0,i}^T \Sigma^{-1} e_i)] $$
Note that, $ \theta_{0,i}^T \Sigma^{-1} e_i = \sum\limits_{ j \neq i } \theta_j t_{j,i} $ and from the independence of $ \theta_{j}$'s we can conclude that, 
$$ E_{ \theta_{0,i}} [ \exp ( - k^2  \theta_{0,i}^T \Sigma^{-1} e_i)] 
= \prod\limits_{j \neq i} E [ e^{ - t_j \theta_j } ] =  \prod\limits_{j \neq i}  ( pe^{ - k^2 t_{j,i} } + (1-p) )     $$ 

Thus, we finally obtain a simplified expression of the optimal test statistic as the following 
$$P( \theta_i = 0 \: | \:  \mathbf{X} ) \: = \:  \frac{1}{1 + \frac{pU_i}{1-p}} $$  

Where $ U_i = \exp ( - ( \frac{k^2}{2} t_{i,i}  -  k \sum\limits_{ j=1}^{n} t_{j,i}  x_j )) ( \prod\limits_{j \neq i}  ( pe^{ - k^2 t_{j,i} } + (1-p) )  )  $  and   $\underset{\sim}{t_i} =(t_{1,i} , ..., t_{n,i} ) $ is the $i$-th column of $\Sigma^{-1}$

\nocite{romano2008control,sarkar2002some,efron2007correlation,efron2010correlated,efron2012large,efron2001empirical,efron2007size}

\bibstyle{plain}
\bibliography{references.bib}

\end{document}